\documentclass[12pt]{amsart}
\usepackage{amscd,amsmath,amssymb,amsfonts}
\theoremstyle{plain}
\newtheorem{thm}{Theorem}
\newtheorem{lem}[thm]{Lemma}
\newtheorem{cor}[thm]{Corollary}
\newtheorem{prop}[thm]{Proposition}

\theoremstyle{definition}

\newtheorem{remark}[thm]{Remark}

\newtheorem{claim}[thm]{Claim}

\numberwithin{thm}{section}
\numberwithin{equation}{section}

\newcommand{\sA}{{\mathcal A}}
\newcommand{\sB}{{\mathcal B}}
\newcommand{\sC}{{\mathcal C}}

\newcommand{\sE}{{\mathcal E}}

\newcommand{\sH}{{\mathcal H}}
\newcommand{\sI}{{\mathcal I}}
\newcommand{\sJ}{{\mathcal J}}
\newcommand{\sK}{{\mathcal K}}

\newcommand{\sO}{{\mathcal O}}

\newcommand{\sT}{{\mathcal T}}

\newcommand{\A}{{\mathbb A}}

\newcommand{\C}{{\mathbb C}}

\newcommand{\E}{{\mathbb E}}

\newcommand{\Q}{{\mathbb Q}}
\newcommand{\R}{{\mathbb R}}

\newcommand{\V}{{\mathbb V}}

\newcommand{\Z}{{\mathbb Z}}
\newcommand{\of}{{\overline{f}}}
\newcommand{\oX}{{\overline{X}}}
\newcommand{\oS}{{\overline{S}}}
\catcode`\@=11
\def\opn#1#2{\def#1{\mathop{\kern0pt\fam0#2}\nolimits}} 
\def\underrightarrow{\mathpalette\underrightarrow@}
\def\underrightarrow@#1#2{\vtop{\ialign{$##$\cr
 \hfil#1#2\hfil\cr\noalign{\nointerlineskip}%
 #1{-}\mkern-6mu\cleaders\hbox{$#1\mkern-2mu{-}\mkern-2mu$}\hfill
 \mkern-6mu{\to}\cr}}}

\def\underleftarrow{\mathpalette\underleftarrow@}
\def\underleftarrow@#1#2{\vtop{\ialign{$##$\cr
 \hfil#1#2\hfil\cr\noalign{\nointerlineskip}#1{\leftarrow}\mkern-6mu
 \cleaders\hbox{$#1\mkern-2mu{-}\mkern-2mu$}\hfill
 \mkern-6mu{-}\cr}}}
\let\amp@rs@nd@\relax
\newdimen\ex@
\newdimen\bigaw@
\newdimen\minaw@
\newdimen\minCDaw@
\newif\ifCD@
\def\minCDarrowwidth#1{\minCDaw@#1}

\def\@CD{\def\A##1A##2A{\llap{$\vcenter{\hbox
 {$\scriptstyle##1$}}$}\Big\uparrow\rlap{$\vcenter{\hbox{$\scriptstyle##2$}}$}&&}%
\def\V##1V##2V{\llap{$\vcenter{\hbox
 {$\scriptstyle##1$}}$}\Big\downarrow\rlap{$\vcenter{\hbox{$\scriptstyle##2$}}$}&&}%
\def\={&\hskip.5em\mathrel
 {\vbox{\hrule width\minCDaw@\vskip3\ex@\hrule width
 \minCDaw@}}\hskip.5em&}%
\def\verteq{\Big\Vert&&}%
\def\noarr{&&}%
\def\vspace##1{\noalign{\vskip##1\relax}}\relax\let\amp@rs@nd@&\iffalse}\fi
 \CD@true\vcenter\bgroup\relax\let\\=\cr\iffalse}\fi\tabskip\z@skip\baselineskip20\ex@
 &\hfill$\m@th##$\hfill\cr}
\def\@endCD{\cr\egroup\egroup}
\def\>#1>#2>{\amp@rs@nd@\setbox\z@\hbox{$\scriptstyle
 \;{#1}\;\;$}\setbox\@ne\hbox{$\scriptstyle\;{#2}\;\;$}\setbox\tw@
 \hbox{$#2$}\ifCD@
 \global\bigaw@\minCDaw@\else\global\bigaw@\minaw@\fi
 \ifdim\wd\z@>\bigaw@\global\bigaw@\wd\z@\fi
 \ifdim\wd\@ne>\bigaw@\global\bigaw@\wd\@ne\fi
 \ifCD@\hskip.5em\fi
 \ifdim\wd\tw@>\z@
 \mathrel{\mathop{\hbox to\bigaw@{\rightarrowfill}}\limits^{#1}_{#2}}\else
 \mathrel{\mathop{\hbox to\bigaw@{\rightarrowfill}}\limits^{#1}}\fi
 \ifCD@\hskip.5em\fi\amp@rs@nd@}
\def\<#1<#2<{\amp@rs@nd@\setbox\z@\hbox{$\scriptstyle
 \;\;{#1}\;$}\setbox\@ne\hbox{$\scriptstyle\;\;{#2}\;$}\setbox\tw@
 \hbox{$#2$}\ifCD@
 \global\bigaw@\minCDaw@\else\global\bigaw@\minaw@\fi
 \ifdim\wd\z@>\bigaw@\global\bigaw@\wd\z@\fi
 \ifdim\wd\@ne>\bigaw@\global\bigaw@\wd\@ne\fi
 \ifCD@\hskip.5em\fi
 \ifdim\wd\tw@>\z@
 \mathrel{\mathop{\hbox to\bigaw@{\leftarrowfill}}\limits^{#1}_{#2}}\else
 \mathrel{\mathop{\hbox to\bigaw@{\leftarrowfill}}\limits^{#1}}\fi
 \ifCD@\hskip.5em\fi\amp@rs@nd@}
\newenvironment{CDS}{\@CDS}{\@endCDS}
\def\@CDS{\def\A##1A##2A{\llap{$\vcenter{\hbox
 {$\scriptstyle##1$}}$}\Big\uparrow\rlap{$\vcenter{\hbox{$\scriptstyle##2$}}$}&}%
\def\V##1V##2V{\llap{$\vcenter{\hbox
 {$\scriptstyle##1$}}$}\Big\downarrow\rlap{$\vcenter{\hbox{$\scriptstyle##2$}}$}&}%
\def\={&\hskip.5em\mathrel
 {\vbox{\hrule width\minCDaw@\vskip3\ex@\hrule width
 \minCDaw@}}\hskip.5em&}
\def\verteq{\Big\Vert&}
\def\novarr{&}
\def\noharr{&&}
\def\SE##1E##2E{\slantedarrow(0,18)(4,-3){##1}{##2}&}
\def\SW##1W##2W{\slantedarrow(24,18)(-4,-3){##1}{##2}&}
\def\NE##1E##2E{\slantedarrow(0,0)(4,3){##1}{##2}&}
\def\NW##1W##2W{\slantedarrow(24,0)(-4,3){##1}{##2}&}
\def\slantedarrow(##1)(##2)##3##4{\thinlines\unitlength1pt\lower 6.5pt\hbox{\begin{picture}(24,18)%
\put(##1){\vector(##2){24}}%
\put(0,8){$\scriptstyle##3$}%
\put(20,8){$\scriptstyle##4$}%
\end{picture}}}
\def\vspace##1{\noalign{\vskip##1\relax}}\relax\let\amp@rs@nd@&\iffalse}\fi
 \CD@true\vcenter\bgroup\relax\let\\=\cr\iffalse}\fi\tabskip\z@skip\baselineskip20\ex@
 &\hfill$\m@th##$\hfill\cr}
\def\@endCDS{\cr\egroup\egroup}
\begin{document}

\title[Weight one vanishing]{Chern classes of Gau{\ss}-Manin bundles of
weight 1 vanish}

\author{H\'el\`ene Esnault}
\address{Universit\"at Essen, FB6, Mathematik, 45117 Essen, Germany}
\email{esnault@uni-essen.de}
\author{Eckart Viehweg}
\address{Universit\"at Essen, FB6, Mathematik, 45117 Essen, Germany}
\email{viehweg@uni-essen.de}

\date{Dec. 31, 2001}
\begin{abstract}
We show that the Chern character of a variation of polarized Hodge
structures of weight one with nilpotent residues at $\infty$
dies up to torsion in the Chow ring, except in codimension 0.
\end{abstract}
\maketitle
\begin{quote}

\end{quote}

\section{Introduction}

Let $f: X\to S$ be a proper smooth family over a smooth base,
defined
over a field $k$ of characteristic 0. Then the Gau{\ss}-Manin
bundles $\sH^i:=R^if_*\Omega^\bullet_{X/S}$ are endowed with the
Gau{\ss}-Manin connection, and thereby, e.g. by  Chern-Weil
theory, their Chern classes die in de~Rham cohomology. By
Griffiths' fundamental theorem \cite{Gr}, the Gau{\ss}-Manin connection is
regular singular, and thinking of $k=\C$, the underlying local
monodromies around the components of $T$ are quasi-unipotent.
If they are unipotent, then Deligne's extension $\overline{\sH}^i$
(\cite{D}) has nilpotent
residues, and an Atiyah class computation (\cite{EV},
appendix B) shows that the Chern classes of $\overline{\sH}^i$
also die in de~Rham cohomology.

On the other hand,
Mumford (\cite{Mu})
remarked that the Grothendieck-Riemann-Roch theorem
applied to the structure sheaf and the dualizing sheaf of
a family of curves $f$,
yields vanishing, up to torsion, of
the Chern classes of $\sH^1$,
a much stronger information
than the vanishing in  de~Rham cohomology.
If one compactifies $f$ as a semistable family
of curves $\bar{f}$, Deligne's extension $\overline{\sH}^1$ is simply
$R^i\bar{f}_*\Omega^\bullet_{X/S}(\log \infty)$ and again,
Grothendieck-Riemann-Roch allows to conclude that the Chern classes
of $\overline{\sH}^1$ die, up to torsion,
in the Chow ring as well (\cite{Mu}). This led the
first author to wonder whether Gau{\ss}-Manin bundles in general
can yield non-trivial algebraic cycles in Chow groups (see
\cite{E}). For example,
it is proven in \cite{BE} that the algebraic Chern-Simons invariants of
Gau{\ss}-Manin bundles in characteristic $p>0$ always die (up to
torsion).

The main theorem of this article is
\begin{thm} \label{mainthm}
Let $B$ be a smooth complex variety, with a compactification
$B\subset \overline{B}$ such that $\overline{B}$ is smooth and
$T:=\overline{B}\setminus B$ is a divisor with normal crossings.
Let $\sH^1$ be a variation of polarized  pure Hodge structures of weight
1, with unipotent local monodromies along the components of $T$,
and let $\overline{\sH}^1$ be its Deligne's extension with
nilpotent residues. Then one has
$$ {\rm ch}(\overline{\sH}^1)\in CH^0(\overline{B})\otimes \Q.$$
\end{thm}

Applying
Grothendieck-Riemann-Roch to powers of a principal polarization
in a family of abelian varieties,
van der Geer proved vanishing for $\sH^1$
(\cite{vdG}). A reduction  to Mumford's curve case,
via the Abel-Jacobi map  for genus $\le 3$  (\cite{vdG2})
and  via the Abel-Prym
map for $g=4,5$ (\cite{I}), yields vanishing for
$\overline{\sH}^1$.

On the other hand, Mumford's argument for a family of abelian varieties shows
immediately that the alternate sum of the Gau{\ss}-Manin bundles
$\sH^i$ has
vanishing Chern classes, up to torsion.
The problem is therefore to find a way to separate
the different weights. Van der Geer's argument seems to be a detour.
De Rham cohomology is not coherent cohomology, yet his proof
relies of the remarkable identity ${\rm Td} \ \E^\vee={\rm ch}
\  (f_*L)$, where $\E^\vee$
is  the Hodge bundle and $L$ is a principal polarization. This identity
does not extend across the boundary, as one easily checks for families of
elliptic curves.

In this note, we present a way to separate the different
weights in the spirit of de~Rham cohomology.
We mod out the family of abelian varieties by the $(-1)$ action,
which has the effect to separate the even from the odd weights. Then
we observe that vanishing for the sum of the Gau{\ss}-Manin
bundles in even weights is equivalent to vanishing for $\sH^1$
(Lemma \ref{calcul}).

The next step is to extend the quotient across
the boundary, keeping track both of the Gau{\ss}-Manin bundles
and of the Riemann-Roch theorem.
To this aim, one has to consider a compactification of the
universal abelian variety with level $n\ge 3$ structure
over the moduli space $\sA_{g,n}$, and to extend the
$(-1)$ action, controlling the fixed points. This is
performed by a careful study of \cite{FC} (see Theorem \ref{thmFC}).
One has then to understand the effect of the boundary on the
Riemann-Roch-Grothendieck theorem. This is Theorem \ref{GRR},
which is perhaps of independent interest. The philosophy of this
log version is that most of the extra terms one has in the
Riemann-Roch formula are killed by the existence of the residues
on the relative log 1-forms.  The other ones die when one assumes
that the relative 1-log forms essentially come from the base of
the family.

A version of Theorem \ref{mainthm} holds over any field of
characteristic $p \neq 2$ (see Remark \ref{kchar0} and
Theorem  \ref{charp}), replacing the Chern character
of the alternate sum of the Gau{\ss}-Manin bundles
by the Chern character of the alternate sum of the cohomology of
the relative differential forms with log poles.

{\it Acknowledgement}: It is a pleasure to thank Luc Illusie for
a discussion on his results in log geometry and Valery Alexeev
for encouraging us to dig out from \cite{FC} some useful
geometric information.

\section{A numerical computation}
Let $\sH$ be a bundle of rank $2g$ over a smooth variety $S$,
which is an extension of a bundle $\E^\vee$ by its dual $\E$.
\begin{lem} \label{calcul} If
\begin{gather}
{\rm ch} (\sum_{i=1}^g \wedge^{2i}\sH)  \in
CH^{0}(S)\otimes \Q, \ {\rm and}\notag  \\
{\rm ch} (\sum_{i=0}^{g-1} \wedge^{2i+1}\sH) =0 \in
CH^{\bullet }(S)\otimes \Q, \notag
\end{gather}
then one has
\begin{gather}
{\rm ch} (\sH)\in CH^{0}(S)\otimes \Q.\notag
\end{gather}
\end{lem}
\begin{proof}
Let $K(S)$ denote the $K$-group of vector bundles on $S$.
Setting as usual (see \cite{H}, \cite{BS}, for example)
\begin{gather}
\lambda_t(\sH)=\sum_{i=0}^{2g}\lambda^i(\sH)t^i \in K(S)[[t]]\
{\rm with} \notag\\
\lambda^i(\sH)=[\wedge^i\sH]\in K(S), \notag
\end{gather}
and  denoting by $e_i,$ $i=1,\ldots, g$, the Chern roots of $\E$,
one has
\begin{gather}
\lambda_t(\sH)=\lambda_t(\E)\cdot \lambda_t(\E^\vee)=\notag
\prod_{i=1}^g (1+e_it)(1+e_i^\vee t).\notag
\end{gather}
Thus
\begin{gather}
{\rm ch}(\lambda_t(\sH))=\prod_{i=1}^g (1+e^{a_i}t)(1+e^{-a_i}t)
\in CH^\bullet(S) \otimes \Q  \notag \\
{\rm with} \ a_i=c_1(e_i)  .\notag
\end{gather}
We set
\begin{gather}
{\rm ch}^{\rm even}(\wedge \sH):=\frac{1}{2} ({\rm ch} \lambda_1(\sH) +
{\rm ch}\lambda_{-1}(\sH)) \notag\\
=\frac{1}{2}\prod_{i=1}^g(1+e^{a_i})(1+e^{-a_i}) +\frac{1}{2}\prod_{i=1}^g
(1-e^{a_i})(1-e^{-a_i}),\notag\\
{\rm ch}^{\rm odd}(\wedge \sH):=\frac{1}{2} ({\rm ch} \lambda_1(\sH) -
{\rm ch}\lambda_{-1}(\sH)) \notag\\
=\frac{1}{2}\prod_{i=1}^g(1+e^{a_i})(1+e^{-a_i}) - \frac{1}{2}\prod_{i=1}^g
(1-e^{a_i})(1-e^{-a_i}).\notag
\end{gather}
Thus the assumption is equivalent to
\begin{gather} \label{sep}
\prod_{i=1}^g(1+e^{a_i})(1+e^{-a_i})=0 \in CH^{\geq 1}(S)\otimes
\Q\notag\\
\prod_{i=1}^g (1-e^{a_i})(1-e^{-a_i})=0 \in CH^\bullet(S) \otimes
\Q.\notag \end{gather}
The first relation  reads
\begin{gather}
\prod_{i=1}^g(1+e^{a_i})^2 e^{-a_i}
(=2^{2g}) \in CH^{0}(S)\otimes
\Q,\notag
\end{gather}
or equivalently
\begin{gather}
-\sum_{i=1}^g a_i + 2 \sum_{i=1}^g \log (1+e^{a_i})
\in CH^0(S)\otimes
\Q.\notag
\end{gather}
Setting
\begin{gather}
\psi(t)=\log (1+e^t) \notag
\end{gather}
one has
\begin{gather}
\psi'(t)=\frac{e^t}{1+e^t}=1 -\varphi(t), \ \text{with} \notag \\
\varphi(t)=\frac{1}{1+e^t}=\frac{1}{2} \sum_{n=0}^\infty E_n(0)
\frac{t^n}{n!},\notag
\end{gather}
where the $E_n(0)$ are the Euler numbers at 0 (see\cite{Ba},
1.14, (2)).

Vanishing of ${\rm ch}(\sH)$ in $CH^{\ge 1}(S)\otimes \Q$ is equivalent to
vanishing of
$${\rm ch}_{2\bullet}(\E) \in CH^{2\bullet}(S)\otimes \Q, \
\bullet\geq 1.$$
Thus it is equivalent to the assertion that none of the
 odd coefficients of the expansion of $\varphi(t)$ is vanishing, that
is that $E_{2n-1}(0)\neq 0$ for all $n\geq 1$. By \cite{Ba}, 1.14,
(7), one has
$$E_{2n-1}(0)=\frac{2(1-2^{2n})}{2n}B_{2n}(0),$$
where $B_n(0)$ are the Bernoulli numbers at 0,  and by
\cite{Ba}, 1.13, (16), $B_{2n}(0)\neq 0$ for all $n\geq 1$.
This concludes the proof.
\end{proof}

\section{The geometry of the compactified family of abelian varieties}
In this section, we extract from \cite{FC} the necessary
geometric information in order to find a model for the
compactification of a family of abelian varieties, which will
allow us to apply in section 5 a log version of the
Grothendieck-Riemann-Roch theorem.

We use the following notations. Fixing the level $n$, we
denote by $S=\sA_{g,n}$ the moduli stack of abelian varieties with
level $n$ structure (\cite{GIT}). For $n\geq 3$, not divisible
by the characteristic of $k$, $S$ is a
scheme and it carries a universal family $f: X\to S$ of abelian
varieties. We consider one of the  compactifications
$\of: \oX\to \oS$  described
on p. 195 of \cite{FC}.

More precisely, the compactification $S\subset \oS$ is
determined by a certain polyhedral cone decomposition
$\{\sigma_\alpha \}$ of $C(N)$, where $N$ is a free abelian
group of rank $g$, $B(N)$ is the space of integer
valued symmetric bilinear forms and $C(N)\subset B(N_{\R})$ is
the convex cone of all positive semi-definite symmetric bilinear
forms whose radicals
$$
{\rm Ker}(b:N_\R \to N_\R^\vee),
$$
are defined over $\Q$ (p. 96, 2.1). The interior
$C^\circ(N)$ consists of the positive definite forms.

Then, $\oS$ being chosen, one considers
$$\tilde{B}(N)=B(N)\times N^\vee.$$
A compactification $\of: \oX\to \oS$
is determined by the choice of a polyhedral cone decomposition
$\{\tau_\beta\}$ of the  cone
\begin{gather*}
\tilde{C}(N)=\{(b,\ell), \ell=0 \ \text {on}  \
{\rm Ker} \ b\}\subset  \  C(N)\times N_{\R}^\vee \\
\subset \tilde{B}(N_{\R})=B(N_\R)\times N_\R^\vee
\end{gather*}
(p. 195, last section).
For the existence of $\of$ one needs (p. 196, 1.3 (v))
that any $\tau_\beta$ maps into a $\sigma_\alpha$. Recall
moreover that one of the conditions on the polyhedral cone
decompositions requires $\{\sigma_\alpha\}$ to be
$GL(N)$-invariant (p. 96, 2.2), $\{\tau_\beta\}$ to be
$GL(N)\ltimes N$-invariant (p. 196, 1.3), and that there are
finitely many orbits.

As is underlined on p. 195, l.1, the family $f: X\to S$ does not
necessarily extend to a semi-abelian group scheme $G \to
\oS$ embedded in $\of: \oX\to
\oS$. Yet, Remark 1.4 p. 197 asserts that it is  possible
to further refine the cone decompositions in such a way that the
natural  section of $\tilde{B}(N)\to B(N)$ respects the cone
decomposition, which means that any $\sigma_\alpha \times \{0\}$
is precisely one of the  $\tau_\beta$. This guarantees that $f:
X\to S$ extends to a semi-abelian group scheme $G\to
\oS$ embedded in $\of: \oX\to
\oS$.

We set $T:= \oS\setminus S, \ Y:=\oX\setminus X$.

Refining, we may assume $\{\tau_\beta\}$ and $\{\sigma_\alpha\}$ to be
smooth (p. 96, 2.3 and p.98, (iii)) and $\{\sigma_\alpha\}$
to satisfy the condition (ii) on p. 97. In particular, this says
that both $\oS, \ \oX$ are smooth, that $T, \  Y$ are
normal crossings divisors and that the components of $T$ are
non-singular (p. 118, 5.8, a).

For $n \geq 3$, it is explained on p.172 and p. 173 how to refine
a given  polyhedral smooth cone decomposition to force
$\oS$ to be a projective and smooth scheme. We remark,
that via p. 173, c), whatever projective polyhedral
decomposition is chosen to define the compactification, it is
always possible to refine it to a finer smooth projective one.
Of course, this changes the compactifications of $S$ and $X$.
For the family $X$, one first quotes Theorem 1.1 p. 195 which
yields $\of:\oX\to \oS$ as a morphism of algebraic stacks.
However, by p. 207, l.4 to 8, we know that $\oX$ is a
smooth projective variety and that $\of$ is then
consequently a projective morphism.

We have essentially reached the first part of the following theorem.
\begin{thm}[Faltings-Chai] \label{thmFC} Let $k$ be a field
containing the $n$-th roots of unity. For $n$ even $ \geq 4$ and
not divisible by the characteristic of $k$, there is a
compactification $\of:\oX\to \oS$ of the universal
family of principally polarized abelian varieties of genus $g$
with level $n$ structure, with the following properties:
\begin{itemize}
\item[1.] $\oX$ and $\oS$ are smooth projective
varieties.
\item[2.]
$T, \ Y$ are normal crossings divisors with smooth irreducible
components.
\item[3.]
The  sheaf
of relative 1-forms $\Omega^1_{\oX/\oS}(\log
Y) $ with logarithmic poles along $Y$ is locally free.
\item[4.] The Hodge bundle
$\E=\of_*\Omega^1_{\oX/\oS}(\log
Y) $ is locally free.
\item[5.]  One has $f^*\E=\Omega^1_{\oX/\oS}(\log
Y) $.
\item[6.]  One has $R^q\of_*(
\Omega^p_{\oX/\oS}(\log Y))=
\wedge^q \E^\vee \otimes \wedge^p\E.$
\item[7.] When the ground field $k$ has characteristic 0,
the  Gau{\ss}-Manin sheaf
$R^q\of_*\Omega^\bullet_{\oX/\oS}(\log
Y)=:\overline{\sH}^q$ is Deligne's extension \cite{D} of its restriction
to $S$. In particular, it is locally free. The residues of the
Gau{\ss}-Manin connection are nilpotent.
\item[8.] $f: X\to S$ extends to a semi-abelian group scheme
$G\to \oS$ embedded into $\of: \oX\to
\oS$.
\item[9.] The level $n$-structure sections $S_i$ of $f: X\to S$ extend
to disjoint sections $\overline{S_i}$ of $ \of: \oX
\to \oS$.
\item[10.] The $\iota:=(-1): X\to X$ involution over $S$ extends to an
involution, still denoted by $\iota: \oX\to
\oX$ over $\oS$.
\item[11.] The fixed points of $\iota$ lie in $\cup
\overline{S_i}$.
\end{itemize}
\end{thm}
\begin{proof}
1., 8. have already been discussed, as well as part of
2. For 3., 4., 5., 6., we
refer to Theorem 1.1, p. 195. For 7., we know (p. 218, (4)) that
the Hodge to de~Rham spectral sequence degenerates, which implies
via 6. that $\overline{\sH}^i$ is locally free. On the other
hand, restricting to a generic curve $\overline{C}$ in $
\oS$ intersecting $T$ in general points $\overline{C}
\setminus C$, one obtains a family $h:\overline{W}\to
\overline{C}$, the fibres of which are all normal crossings
divisors. As well known (see \cite{Gri}, p. 130, for example)
$\overline{\sH}^i$ is the Deligne extension of
$\overline{\sH}^i|_C$. On the other hand, 8. implies that the
Gau{\ss}-Manin connection on $\overline{\sH}^i|_C$ has nilpotent
residues.

For 9., 10., 11. and to see that the components of $Y$ can be
assumed to be smooth, we need a more precise discussion of the
polyhedral cone decomposition and its relation with the
compactification. Obviously the three properties hold true over
$S$, hence extending them to the boundary is a local question. As
on p. 207, 2., we can even replace $\oS$ by the formal
completion along a certain stratum $\mathfrak{Z}$ and
$\of:\oX\to \oS$ by the pullback
family. Doing so, we are allowed to use the description of
$\of$ given  p. 201 and p. 203.

Recall that the toroidal embeddings $\overline{F}\to
\overline{E}$ given by the polyhedral cone decompositions are
stratified by locally closed subschemes $\mathfrak{Z}_{\tau_\beta}$ and
$\mathfrak{Z}_{\sigma_\alpha}$ (p. 100, 2.5, (iv)). The $\mathfrak{Z}_{\sigma_\alpha}$
and $\mathfrak{Z}_{\tau_\beta}$ are orbits under the torus action, and their
codimension is equal to the dimension of the $\R$-vectorspace
spanned by $\sigma_\alpha$ or $\tau_\beta$, respectively.

If the fibre $G_0$ of $G$ over the general point of the stratum
$\mathfrak{Z}$ is a torus, then the formal completions along the stratum,
together with the pullback of $\of:\oX\to
\oS$, are obtained by restricting $\overline{F}$ to the
formal completion of $\overline{E}$ along a stratum
$\mathfrak{Z}_{\sigma_\alpha}$, with $\sigma_\alpha \subset C^\circ(N)$, and
by taking the quotient by $N$ (p. 201, last section). In general
$G_0$ will be an extension of an abelian variety
$A$ of dimension $g-r$ by a torus. As sketched on p. 202 - 203, one has to
replace $N$ by an $r$-dimensional quotient lattice $N_{\xi}$ in
this case. In particular one may again assume that
$\sigma_\alpha$ lies in the interior of the cone $C(N_\xi)$.
Since the combinatorial description remains the same, we drop the
${}_\xi$.

By p. 100, 2.5 (i) and (ix), the category of rational partial
polyhedral cone decompositions is equivalent to the category of
torus embeddings. Thus composing a given torus embedding with
$\iota=(-1)$ corresponds to the action of $\iota$ on the data
giving the polyhedral cone decomposition. So we just have to
verify that $\iota$ respect those data. $\iota$ acts on $N$ by
multiplication with $(-1)$, hence the action is trivial on
$B(N)$, and is $(-1)$ on $N^\vee$. An element $\mu \in N$ acts on
$(b,\ell) \in \tilde{B}(N)=B(N)\times N$ by
$$
(b,\ell) \mapsto (b,\ell + b(\mu, \ \ )),
$$
while $\gamma \in GL(N)$ acts via
$$ (b,\ell) \mapsto (b\circ (\gamma^{-1}, \gamma^{-1}), \ell
\circ \gamma^{-1}) $$
(p. 196, first section). In particular,
$$ (-{\rm Id}, 0)  \in GL(N) \ltimes N$$
maps
$$(b, \ell)\in B(N)\times N^\vee$$
to $(b,-\ell).$ Since the polyhedral cone decompositions are
invariant under $GL(N) \ltimes N$, it is invariant under $\iota$.

As for 9., we just remark that a morphism from $\oS$ to
the corresponding compactification $\oS_1$ of the
moduli stack $\sA_{g,1}$ is defined by multiplication with $n$ on
the torus (p. 130, 6.7. (6)). In different terms, one keeps the
cone decompositions in $B(N_\R)$ and
$\tilde{B}(N_\R)\times N_\R$, but changes the integral
structure by multiplication with $n$ on $N$. To see that the
closure of the sections of $f:X\to S$ of order $n$ in $\oX$
are disjoint sections of $\of$, it is again sufficient
to consider the pullback of $\of$ to formal completions
of the strata in $\oS$. By p. 202, first section,
the $n$-torsion points are the pull-back of the zero-section
of the semi-abelian group scheme over the formal completion.

As already seen, $\iota$ acts trivially on the cone
$\sigma_\alpha\subset B(N_\R),$
and the fixed points under the $\iota$ involution lie in strata
$\mathfrak{Z}_{\tau_\beta}$ of $\overline{F}$, the $N$-orbit of which are invariant
under $\iota$. We assumed that $\{\tau_\beta\}$ is smooth, that
is each cone $\tau_\beta$ is generated by a partial $\Z$-basis
$$\big((b_1,\ell_1),\ldots, (b_r,\ell_r)\big) \text{ \ of \ }
B(N)\times N^\vee.$$
Thus one has a $\mu \in N$ with
$$(b_i,-\ell_i)=(b_{j(i)}, \ell_{j(i)} +b_{j(i)}(\mu, \ \ ))$$
for $i=1,\ldots, r$. We have taken an even level $n$. Thus
$N=2\cdot N'$ for another integral structure, and one has
$\mu=2\cdot \mu'$ for a $\mu'\in N'$. We obtain
$$b_i=b_{j(i)}\text{ \  and \ }\ell_i +\ell_{j(i)}=2 b_{j(i)}(\mu',
\ \ ).$$
Twisting the free $\Z$-module $B(N)\times N^\vee$ by
$\Z/2$ over $\Z$, the basis elements $(b_i, \ell_i)$ and
$(b_{j(i)}, \ell_{j(i)})$ become equal, which implies that
$j(i)=i$. This in turn implies that $2\ell_i=2b_i(\mu', \ \ )$,
and shows that the dimension of the subspace of $B(N)\times
N^\vee$ generated by $\tau_\beta$ is the same as the dimension
of its image in $B(N)$. This finally implies that ${\rm codim}
(\mathfrak{Z}_{\tau_\beta})$ is equal to the codimension of its image in
$\overline{E}$. So the fixed points of $\iota$ all lie in the
smooth locus of $\of$.

We verified all the conditions stated in \ref{thmFC}, except that
$Y$ can still have singular components. However, since 1., 3. -
11. are compatible with the blowing up of non-singular strata
of the singular locus of $Y_{\rm red}$, this last point can
be achieved.
\end{proof}

\section{A log version of the Grothendieck-Riemann-Roch theorem}
In this section, we show that the Grothendieck-Riemann-Roch
theorem extends to a log version.
\begin{thm} \label{GRR} Let $\of: \oX \to
\oS$ be a projective morphism of relative dimension $g$ over a
field, with $\oX, \ \oS$ smooth, compactifying the smooth
projective morphism $f: X\to
S$, with the
following properties:
\begin{itemize}
\item[1.]
Both $T=\oS\setminus S$ and $Y:=(\of^{*}(T))_{\rm
red}$ are normal crossings divisors with smooth irreducible components.
\item[2.] The sheaves $\Omega^i_{\oX/\oS}(\log Y)$ are
locally free.
\item[3.]
There are  cycles $W\in CH^g(\oX)\otimes \Q, \ \xi \in
CH^g(\oS)\otimes \Q$, such that
$$ c_g(\Omega^1_{\oX/\oS}(\log Y))= \of^*(\xi) + W \in
CH^g(\oX)\otimes \Q,$$
with the property $ W\cdot Y_i=0 \in CH^{g+1}(\oX) \otimes \Q$
for all irreducible components of $Z:=\of^*(T) - Y$.
\end{itemize}
Then one has
\begin{gather}
\sum_i (-1)^i {\rm ch}(\sum_{j}
R^j\of_*\Omega^{i-j}_{\oX/\oS}(\log Y)) \in CH^0(\oS)\otimes
\Q.\notag
\end{gather}
\end{thm}
\begin{remark} \label{rmk:GRR}
Theorem \ref{GRR} applies in particular when
$Z=\emptyset$, that is when the fibers have no multiplicities.
In this case, 3. is automatically fulfilled and, as we will see,
the proof does not require any combinatorics.
\end{remark}
\begin{proof}
We apply the Grothendieck-Riemann-Roch theorem \cite{BS}
to the alternate sum of the sheaves
$\Omega^i_{\oX/\oS}(\log Y)$.  It yields
\begin{gather}
\sum_i (-1)^i {\rm ch}
(\sum_{j}
R^j\of_*\Omega^{i-j}_{\oX/\oS}(\log Y))\notag \\
=\of_* {\rm
Todd}(T_{\oX/\oS})\cdot (\sum_i (-1)^i
{\rm ch} (\Omega^i_{\oX/\oS}(\log Y))).\notag
\end{gather}
We consider the residue sequences
\begin{gather}
0 \to {\of}^*\Omega^1_{\oS} \to {\of}^*\Omega^1_{\oS}(\log T) \to
\sO_{\tilde{T}} \to 0 \notag \\
0 \to \Omega^1_{\oX} \to \Omega^1_{\oX}(\log Y) \to \sO_{\tilde Y}
\to 0 ,\notag
\end{gather}
where, in order to simplify notations, we have set
$$
\tilde{T_j} = \of^*T_j \text{ \ \ and \ \ }
\sO_{\tilde{T}}=\bigoplus_j \sO_{\tilde{T}_j},
$$
for the irreducible components $T_j$ of $T$, and similarly
$\sO_{\tilde{Y}}=\oplus_i \sO_{Y_i}$.
Multiplicativity of the Todd class implies
\begin{gather}
{\rm Todd} (T_\oX) \cdot
{\rm Todd} (\of^*T_\oS)^{-1}= {\rm Todd}
(T_{\oX/\oS}(\log Y)) \cdot \notag \\
\prod_{q\ge 1}{\rm Todd} (\sE xt^q(\sO_{\tilde{Y}},
\sO_\oX))^{(-1)^{q+1}} \cdot \prod_{q\ge 1}
{\rm Todd} (\sE xt^q(\sO_{\tilde{T}}, \sO_\oX))^{(-1)^q}.
\notag \end{gather}
We define $\sB$ and $\sC$ via the exact sequences
\begin{gather}
0\to \sB \to \sO_{\tilde{Y}} \to \sO_{\tilde{Y}}/\sO_{\tilde{T}} \to
0\notag \\
0\to \sC \to \sO_{\tilde{T}} \to \sB\to 0.\notag
\end{gather}
Using again multiplicativity, one obtains
\begin{gather}
{\rm Todd} (T_\oX) \cdot {\rm Todd} (\of^*T_\oS)^{-1}= {\rm
Todd}(T_{\oX/\oS}(\log Y) )\cdot \notag \\
 \prod_{q\ge 1}
{\rm Todd} (\sE xt^q(\sO_{\tilde{Y}}/\sO_{\tilde{T}},
\sO_\oX))^{(-1)^{q+1}} \cdot \prod_{q\ge 1}
{\rm Todd} (\sE xt^q(\sC, \sO_\oX))^{(-1)^q}.
\notag \end{gather}
On the other hand, one has the well known  relation \cite{H}
\begin{gather}
{\rm Todd}(T_{\oX/\oS}(\log Y) )\cdot
(\sum (-1)^i{\rm ch} (\Omega^i_{\oX/\oS}(\log Y)))= \notag \\
(-1)^g c_g(\Omega^1_{\oX/\oS}(\log Y) ).\notag
\end{gather}
Let us write
\begin{gather}
\prod_{q\ge 1}
{\rm Todd} (\sE xt^q(\sO_{\tilde{Y}}/\sO_{\tilde{T}},
\sO_\oX))^{(-1)^{q+1}}= 1 + {\rm
err}\notag \\
{\rm with} \ \ {\rm err} \in
CH^\bullet(\oX).\notag
\end{gather}
The sheaf $\sO_{\tilde{Y}}/\sO_{\tilde{T}}$
is a direct sum of structure sheaves $\sO_{Y_{i_1}\cap \ldots
\cap Y_{i_\ell}}$.
The surjection $\Omega^1_{\oX/\oS}(\log Y) \twoheadrightarrow
\sO_{Y_{i_1}\cap \ldots
\cap Y_{i_\ell}}$ implies that
$$c_g(\Omega^1_{\oX/\oS}(\log Y))\cdot( Y_{i_1}\cap \ldots
\cap Y_{i_\ell})=0\in CH^{g+\ell}(\oX).$$
Since ${\rm err}$ is a sum of terms supported in those strata
$Y_{i_1}\cap \ldots
\cap Y_{i_\ell}$, we conclude
\begin{gather}
c_g(\Omega^1_{\oX/\oS}(\log Y)) \cdot
\prod_{q\ge 1}
{\rm Todd} (\sE xt^q(\sO_{\tilde{Y}}/\sO_{\tilde{T}},
\sO_\oX))^{(-1)^{q+1}} \notag\\
=
c_g(\Omega^1_{\oX/\oS}(\log Y)).\notag
\end{gather}
Thus
\begin{gather}
\of_* {\rm
Todd}(T_{\oX/\oS})\cdot (\sum_i (-1)^i
{\rm ch} (\Omega^i_{\oX/\oS}(\log Y)))=\label{form1}\\
(-1)^g\of_*(c_g(\Omega^1_{\oX/\oS}(\log Y))\cdot
\prod_{q\ge 1}{\rm Todd} (\sE xt^q(\sC, \sO_\oX))^{(-1)^q}). \notag
\end{gather}
If $Z=\emptyset$, the sheaf $\sC$ is zero,
hence ${\rm Todd}(\sE xt^q(\sC, \sO_\oX))=1$.
This finishes the proof of Remark \ref{rmk:GRR}.

For the general case we have to study the residue map
and the sheaf $\sO_{\tilde{Y}}/\sO_{\tilde{T}}$
more precisely. Let us first fix notations.
We write
$$ Y=\sum_{i\in \sI} Y_i, \ \ \ T=\sum_{j\in \sJ} T_j,
$$
where the $Y_i, \ T_j$ are prime divisors.
For $I\subset \sI, \  J\subset \sJ$, we define
$$
Y_I=\bigcap_{i\in I} Y_i, \ \ \ T_J=\bigcap_{j\in J} T_j.
$$
Since $\Omega^1_{\oX/\oS}(\log Y)$ is locally
free, $\of$ sends strata to strata. Indeed, choose $J$ maximal
with the property that $T_J \supset \of(Y_I)$. If $Y_I\to T_J$
is not surjective, we find a local parameter $t$ on $\oS$, in
a general point of $\of(Y_I)$, such that $\of(Y_I)$ lies in the
zero locus of $t$ but not in $T_J$. As $\sO_\oS\cdot dt$ is
locally split in $\Omega^1_\oS(\log T)$, and the injection
$$\of^*\Omega^1_\oS(\log T)\subset \Omega^1_{\oX}(\log Y)$$
is locally split, $\sO_\oX \cdot dt$ is locally split in
$\Omega^1_{\oX}(\log Y)$ as well. But since by assumption,
$$t=\sum_{i\in I} y_i \alpha_i, \ \alpha_i\in \sO_{\oX},$$
one finds
$$dt=\sum_{i\in I} y_i\cdot (d\alpha_i +\frac{dy_i}{y_i}
\alpha_i),$$
a contradiction.

This allows to define, for each $I\subset \sI$, the index set
$J:=J(I)\subset \sJ$ with the property
$\of(Y_I)= T_J$.
It yields  a $\Q$-linear map
\begin{gather}
\varphi_I: \bigoplus_{j\in J} \Q\cdot T_j \to \bigoplus_{i\in I}
\Q\cdot Y_i \notag \\
\varphi_I(T_j)=\sum_{i\in I} \nu_i^j Y_i \ \text{ \ where \ } \
\of^*T_j=\sum_{i\in \sI} \nu_i^j Y_i. \notag
\end{gather}
One defines
\begin{gather}\label{defn:delta}
\delta_Y={\rm codim}({\rm Im}(\varphi_I))=|I|- {\rm dim}({\rm
Im}(\varphi_I)).
\end{gather}
\begin{claim} \label{lem2} One has $\delta_I >0$ if and only of
the residue map $$\Omega^1_X(\log Y) \twoheadrightarrow \sO_{Y_I}$$ factors
through $\Omega^1_{\oX/\oS}(\log Y)\twoheadrightarrow \sO_{Y_I}$.
\end{claim}
\begin{proof}
We fix $I, J=J(I)$ as before. The question being local at
generic points of $T_J$ and $Y_I$, we choose local parameters
$\{t_j\}$  and $\{y_i\}$ defining $T_j$ and $Y_i$ locally. Then the
equation of the morphism $\of$ is simply $t_j=\prod_{i\in I}
y_i^{\nu_i^j}$.

As in the beginning of the proof of Theorem \ref{GRR},
functoriality maps
$$
\of^*\Omega^1_\oS(\log T) \xrightarrow{{\rm
res}} \of^*(\bigoplus_{j\in J} \sO_{T_j})
\text{ \ \ to \ \ }
\Omega^1_\oX(\log Y) \xrightarrow{{\rm
res}} \bigoplus_{i\in I} \sO_{Y_i}.
$$
It induces
\begin{gather}
\bigoplus_{j\in J}\of^*(\sO_{T_j})|_{Y_I}=
\bigoplus_{ J}\sO_{Y_I} \to \bigoplus_{I}\sO_{Y_I}\hspace{2cm}\notag \\
\hspace{2.3cm}{\rm res} (\frac{dt_j}{t_j}) \mapsto \sum_{i\in I}
\nu_i^j {\rm res}(\frac{dy_i}{y_i}). \notag
\end{gather}
\end{proof}
\begin{claim} \label{cor3}
Let $i_1,\ldots, i_\mu \in \sI$, not necessarily pairwise
distinct. Assume that $Y_I\neq \emptyset$, and $\delta_I >0$,
where $I\subset \sI$ is the smallest subset of $\sI$ containing
$\{i_1,\ldots , i_\mu\}$. Then the cycle $Y_{i_1}\cdots Y_{i_\mu}
\in CH^\mu(\oX)$ fulfills
\begin{gather}
c_g(\Omega^1_{\oX/\oS}(\log Y))\cdot Y_{i_1} \cdots Y_{i_\mu} =0
\in CH^{g+\mu}(\oX). \notag
\end{gather}
\end{claim}
\begin{proof}
Without loss of generality, we may assume that $I=\{i_1,\ldots,
i_\eta\}$ are pairwise distinct. By \ref{lem2}, one has
$c_g(\Omega^1_{\oX/\oS}(\log Y))\cdot Y_I=0$, thus a fortiori
$c_g(\Omega^1_{\oX/\oS}(\log Y))\cdot Y_I\cdot Y_{i_{\eta
+1}}\cdots Y_{i_\mu} =0$.
\end{proof}
Recall that for the proof of Theorem \ref{GRR}, it remains to show that
the term
$$
(-1)^g\of_*(c_g(\Omega^1_{\oX/\oS}(\log Y))\cdot
\prod_{q\ge 1}{\rm Todd} (\sE xt^q(\sC, \sO_\oX))^{(-1)^q})
$$
in formula (\ref{form1}) is equal to
$(-1)^g\of_*(c_g(\Omega^1_{\oX/\oS}(\log Y)))$.
Writing
$$\tilde{T}_j = \of^*(T_j) \text{ \ \ and \ \
}Z_j=\tilde{T}_j-(\tilde{T}_j)_{{\rm red}}$$
one has a commutative diagram
$$
\begin{CDS}
0\>>> \bigoplus_j\sO_{Z_j}(-(\tilde{T}_j)_{\rm red}) \>>>
\bigoplus_j \sO_{\tilde{T}_j} \>>> \bigoplus_j
\sO_{(\tilde{T}_j)_{ {\rm red}}}
\>>> 0\\
\noharr \V\rho_1 VV \novarr \V = VV \novarr \V V \rho_2V\\
0\>>> \sC \>>> \bigoplus_j \sO_{\tilde{T}_j} \>>>
\bigoplus_i\sO_{Y_i}
\end{CDS}
$$
with exact rows. Since $\rho_1$ is injective ${\rm Ker} \ \rho_2
= {\rm Coker} \ \rho_1=:\sK.$ Thus
\begin{gather}
\prod_{q\ge 1}{\rm Todd}(\sE xt^q(\sC, \sO_\oX))^{(-1)^q}=\notag\\
\prod_{q\ge 1} {\rm Todd} (\sE xt^q(\oplus_j
\sO_{Z_j}(-(\tilde{T}_j)_{{\rm red}}, \sO_\oX)))^{(-1)^q}\cdot
\prod_{q\ge 1}{\rm Todd}(\sE xt^q(\sK, \sO_\oX))^{(-1)^q}.\notag
\end{gather}
One has
\begin{gather}
\sE xt^1(\sO_{Z_j}(-(\tilde{T}_j)_{{\rm red}}),\sO_\oX) =
\sO_\oX(\tilde{T}_j)/\sO_\oX((\tilde{T}_j)_{ {\rm
red}}),\notag \\
\text{while \ \ } \
 \sE xt^q(\sO_{Z_j}(-(\tilde{T}_j)_{{\rm red}}),\sO_\oX)=0, \
q\ge 2 .\notag
\end{gather}
Thus
\begin{gather}
\prod_{q\ge 1} {\rm Todd} (\sE xt^q(
\sO_{Z_j}(-(\tilde{T}_j)_{{\rm red}}, \sO_\oX)))^{(-1)^q}=
1 + v_j,\notag
\end{gather}
where $v_j\in CH^{\ge 1}(\oX)$
can be written as a sum of terms
$$m_{i_1,\ldots,
i_\mu}Y_{i_1}\cdot\cdots\cdot Y_{i_\mu},$$
where at least one of
the $Y_{i_p}$ lies in $Z_j$, and $m_{i_1,\ldots, i_\mu} \in \Z$.
On the other hand, one has the exact sequence
\begin{gather}
0\to \sK\to \bigoplus_j \sO_{(\tilde{T}_j)_{{\rm red}}} \to
\sO_{\tilde{Y}} \to \sO_{\tilde{Y}}/\sO_{\tilde{T}} \to 0.\notag
\end{gather}
Let us write
\begin{gather}
Y=\Phi + V \text{ \ \ and \ \ }
(\tilde{T}_j)_{{\rm red}}=\Phi_j + V_j,\notag
\end{gather}
where each component of $\Phi_j$ maps surjectively onto $T_j$,
where each component of $V$ and $V_j$ maps to a lower
dimensional strata in $T$, and with $\Phi=(\sum_j\Phi_j)_{\rm
red}$. Of course $V_j \subset V$ and $V \subset Z$. Denoting as
usual by $\tilde{ \ }$ the normalization, one has
\begin{gather}
\prod_{q\ge 1}{\rm Todd}(\sE xt^q(\sK, \sO_\oX))^{(-1)^q}=\notag
\\
\prod_{q\ge 1} \big[{\rm Todd}(\sE xt^q(\oplus_j (\sO_{\tilde{\Phi}_j}
\oplus \sO_{\tilde{V}_j}), \sO_\oX))^{(-1)^q}\notag \\
\cdot {\rm Todd}(\sE xt^q( \oplus_j ((\sO_{\tilde{\Phi}_j}
\oplus \sO_{\tilde{V}_j})/\sO_{(\tilde{T}_j)_{\rm red}}, \sO_\oX)))^{(-1)^{q+1}}
\notag \\
\cdot {\rm Todd}(\sE xt^q(\sO_{\tilde{Y}},
\sO_\oX))^{(-1)^{q+1}}
\cdot {\rm Todd}(\sE xt^q(\sO_{\tilde{Y}}/\sO_{\tilde{T}} ,
\sO_\oX))^{(-1)^q}\big]. \notag
\end{gather}
Next, for $i \in \sI$ let us define
$$
N(i) = |\{j\subset \sJ; \ Y_i \subset (\tilde{T}_j)_{\rm red}\}| -1.
$$
Remark that $N(i)$ is zero except when  $Y_i \subset V$.
One has
\begin{gather}
{\rm Todd}(\sE xt^q(\oplus_j (\sO_{\tilde{\Phi}_j}
\oplus \sO_{\tilde{V}_j}), \sO_\oX))^{(-1)^q}\cdot {\rm Todd}(\sE
xt^q(\sO_{\tilde{Y}},\sO_\oX))^{(-1)^{q+1}}=\notag \\
{\rm Todd} (\sE xt^q(\oplus_{i\in \sI} \sO_{Y_i}^{\oplus N(i)},
\sO_\oX))^{(-1)^q}.\notag \end{gather}
Moreover, one has
\begin{gather}
{\rm Todd}(\sE xt^q((\sO_{\tilde{\Phi}_j}
\oplus \sO_{\tilde{V}_j})/\sO_{(\tilde{T}_j)_{\rm red}},
\sO_\oX))^{(-1)^{q+1}} = 1+ w_j,\notag
\end{gather}
where $w_j\in CH^{\ge 1}(\oX)$ is the
sum of terms $m_{i_1,\ldots, i_\mu}Y_{i_1}\cdot\cdots\cdot Y_{i_\mu}$,
where at least two of the $i_p$ are different.

Assume that
$Y_{i_1},\ldots,Y_{i_\mu}$ are all components of $\Phi_j$. The image of
$Y_{\{i_1,\ldots,i_\mu\}}$ is one of the strata of $T$. Since $Y$ is a
normal crossing divisor, $\of (Y_{\{i_1,\ldots,i_\mu\}}) = T_j$.
So $\delta_{\{i_1,\ldots, i_\mu\}}>0$ and by Claim 
\ref{cor3} the intersection of
$c_g(\Omega^1_{\oX/\oS}(\log Y))$ with the cycle
$Y_{i_1}\cdot\cdots\cdot Y_{i_\mu}$ is zero.

Altogether, one has
\begin{gather}
c_g(\Omega^1_{\oX/\oS}(\log Y))\cdot
\prod_{q\ge 1}{\rm Todd} (\sE xt^q(\sC, \sO_\oX))^{(-1)^q}=\notag\\
c_g(\Omega^1_{\oX/\oS}(\log Y))\cdot
(1 + w),\notag
\end{gather}
where the cycle $w \in CH^{\geq 1}(\oX)$ lies in the subspace generated by
products $Y_{\iota_1}\cdot\cdots\cdot Y_{\iota_\ell}$ with
at least one $Y_{\iota_\nu}$ contained in $Z$.

We conclude the proof of Theorem \ref{GRR} applying Proposition \ref{prop:van}.
\end{proof}
\begin{prop} \label{prop:van}.
Let $\of: \oX \to \oS$ be as in Theorem \ref{GRR}. $\ell
\geq 1$ irreducible components $Y_{\iota_1},\ldots,Y_{\iota_\ell}$
of $Y$, with $Y_{\iota_\ell} \subset |Z|$, fulfill
\begin{gather}
\of_*(c_g(\Omega^1_{\oX/\oS}(\log Y)) \cdot
Y_{\iota_1}\cdot\cdots\cdot Y_{\iota_\ell})=0 \in
CH^\ell(\oS)\otimes \Q. \notag \end{gather}
\end{prop}
By Claim \ref{cor3}, we only have to study  strata for which
$\delta_I=0$.
\begin{lem} \label{lem4}
If $\delta_I=0$ and $i\in I$ is given, then there is a divisor
$\Gamma$ supported in $T\subset \oS$, there are indices
$\ell_1,\ldots, \ell_m \in \sI\setminus I$ and multiplicities
$\beta_1,\ldots, \beta_m \in \Q$ fulfilling
\begin{gather}
Y_i\cdot Y_I = \of^*(\Gamma)\cdot Y_I +\sum_{\nu=1}^m \beta_\nu
Y_{I\cup\{\ell_\nu\}} \in CH^{1+|I|}(\oX)\otimes \Q.\notag
\end{gather}
\end{lem}
\begin{proof}
By formula (\ref{defn:delta}), we know that $Y_i$ lies in
the image of $\varphi_I$. Thus  there is a $\Q$-divisor $\Gamma$
supported on $T_J$ such that $Y_i -\of^*(\Gamma)$
is supported away of $Y_a, a\in I$.
\end{proof}
\begin{lem} \label{lem5}
For $\ell \geq 1$ let $\mathfrak{P}_\ell$ be the set of all sets
of indices $I$ with $|I|=\mu\le \ell$,  with $Y_I \neq
\emptyset$, $\delta_I=0$, and such that there is one $i\in I$
with $Y_i\subset |Z|$. For irreducible components
$Y_{\iota_1},\ldots,Y_{\iota_\ell}$ of $Y$, with $Y_{\iota_\ell}
\subset |Z|$, there exist multiplicities $a_I\in \Q$ for $I\in
\mathfrak{P}_\ell$, cycles
$$
\Gamma_I \in CH^{\ell-|I|} (\oS)\otimes \Q,\text{ \ \ with \ \
} \Gamma_I|_{ \oS\setminus T}=0 \in CH^{\ell-|I|}(\oS\setminus
T)\otimes \Q,
$$
such that
\begin{gather} \label{eqn1}
 c_g(\Omega^1_{\oX/\oS}(\log Y))\cdot
Y_{\iota_1}\cdot\cdots\cdot Y_{\iota_\ell}=
c_g(\Omega^1_{\oX/\oS}(\log Y))\cdot
\sum_{I\in \mathfrak{P}_\ell} a_I \of^*(\Gamma_I)\cdot Y_I.
\end{gather}
\end{lem}
\begin{proof}
Obviously $\delta_{\{i\}}=0$, thus $\mathfrak{P}_1=\{i \in \sI,
Y_i\subset |Z|\}$. So for $\ell=1$ there is nothing to show.

We assume now by induction on $\ell$ that formula (\ref{eqn1})
holds true for $\ell$. We want to prove it for $\ell+1$.
It remains to show that for $i\in \sI$ and $I\in
\mathfrak{P}_\ell$,
\begin{gather} \label{eqn3}
c_g(\Omega^1_{\oX/\oS}(\log Y))\cdot Y_i\cdot Y_I
\end{gather}
has the shape required on the right hand side of formula
(\ref{eqn1}).

If $i\notin I$, then either
$\delta_{ \{i \}\cup I}=0$, which implies $\{i\}\cup I \in
\mathfrak{P}_{\ell+1}$, or else the expression (\ref{eqn3})  
is zero by Claim \ref{cor3}.

Thus we assume that $i\in I$.  By Lemma \ref{lem4}, one has
\begin{gather}
Y_i\cdot Y_I=\of^*(\Gamma)\cdot Y_I +\sum_{\nu=1}^m \beta_{\nu}
Y_{\{\ell_\nu\} \cup I},\notag
\end{gather}
with $\ell_\nu \notin I$. Thus again, either $\{\ell_\nu\}\cup
I\in \mathfrak{P}_{\ell+1}$, in which case we are done, or else
$$c_g(\Omega^1_{\oX/\oS} (\log Y))\cdot Y_{\{\ell_\nu\}\cup
I}=0.$$
\end{proof}
\begin{proof}[Proof of Proposition \ref{prop:van}.]
The definition of $\mathfrak{P}_\ell$ implies that the restriction
to $\oX\setminus |Z|$ of the cycles $Y_I\in CH^{|I|}(\oX)\otimes
\Q$ die. In particular one has
$$
W\cdot Y_I=0\in CH^{g+|I|}(\oX)\otimes \Q.
$$
Thus formula (\ref{eqn1}) implies
\begin{gather} \label{eqn4}
c_g(\Omega^1_{\oX/\oS}(\log Y))\cdot Y^\ell=\\
\of^*(\xi)\cdot\sum_{I\in \mathfrak{P}_\ell}
a_I \of^*(\Gamma_I)\cdot Y_I \in CH^{g+\ell}(\oX) \otimes \Q.\notag
\end{gather}
On the other hand, $\delta_I=0$ implies that $|I|\le |J(I)|$.
This in turn implies that the fiber dimension of $Y_I\twoheadrightarrow
T_{J(I)}$ is at least $g$, thus by projection formula, $\of_*$
of the right hand side of formula (\ref{eqn4}) vanishes as well.
This concludes the proof of the proposition.
\end{proof}

Theorem \ref{GRR} implies
\begin{cor} \label{char0}
The assumptions are as in Theorem \ref{GRR}, and
moreover the following conditions are fulfilled:
\begin{itemize}
\item[5.] The characteristic of the ground field $k$ is 0.
\item[6.]
The Hodge to de~Rham spectral sequence $$E_2^{pq}=R^q\of_*
\Omega^p_{\oX/\oS}(\log Y) \Longrightarrow
\overline{\sH}^i:=R^i\of_*\Omega^\bullet_{\oX/\oS}
(\log Y)
$$ degenerates in $E_2$ and the Gau{\ss}-Manin sheaves $\overline{\sH}^i$
are locally free.
\end{itemize}
Then 
\begin{gather}
\sum_i (-1)^i {\rm ch}(\overline{\sH}^i) \in CH^0(\oS)\otimes
\Q.\notag \end{gather}
\end{cor}

\section{Vanishing of the Chern character of the Gau{\ss}-Manin
bundles of weight one.}
In this section, we prove Theorem \ref{mainthm}, which is the
main result of this article.
\begin{remark} \label{kchar0}
Any variation of polarized pure Hodge structures of weight 1 is the
Gau{\ss}-Manin bundle of a smooth polarized family $f: X\to S$ of
abelian varieties (\cite{D2}), with $B=S$. Thus Theorem \ref{mainthm}
may be reformulated over any field of characteristic 0 by saying
that
$$ {\rm ch}(\overline{\sH}^1) \in CH^0(\overline{B})\otimes \Q,$$
where $\sH^1=R^1f_*\Omega^\bullet_{X/S}$ is the first
Gau{\ss}-Manin bundle, and $\overline{\sH}^1$ its Deligne extension with
nilpotent residues.
\end{remark}
\begin{proof}[Proof of Theorem \ref{mainthm}.]
Let $B \subset \overline{B}$ be as in the theorem. There exists a
generically finite morphism $\pi: B_1\to \overline{B}$ together
with a morphism $\psi: B_1\to \oS$, where $\oS$ is as in Theorem
\ref{thmFC}, with
$$\pi^*\sH^1= \psi^* R^1\of_*(\Omega^\bullet_{\oX/\oS}(\log
Y)).$$
Indeed, the polarization of $\sH^1$ is not necessarily coming
from a principal polarization on the underlying family of
abelian varieties, and also, this family might not have a level
$n$-structure, but both can be achieved after replacing $B$ by a
generically  finite covering.

Since $\oS$ and $\overline{B}$ are smooth, projection formula for
$\psi$ and $\pi$ implies that
$$CH(\oS)^\bullet\otimes \Q \text{ \ \ and \ \ }
CH(\overline{B})^\bullet\otimes \Q$$
are subrings of $CH^\bullet(B_1)\otimes \Q$, thus we may assume
that $B= \oS$ and
$$\overline{\sH}^1=R^1\of_*(\Omega^\bullet_{\oX/\oS}(\log Y)).$$
By Corollary \ref{char0}, we know
$$
\sum_i (-1)^i {\rm ch}(\overline{\sH}^i) \in CH^0(\oS)\otimes
\Q.
$$
On the other hand, consider as in Theorem \ref{thmFC},
11., the involution $\iota$ on $\oX$. By loc. cit. 9., 11., the
fixed points of $\iota$ lie in disjoint sections
$\overline{S_\alpha}\subset \oX$ of $\of$. In particular, setting again
$Y=(\of^*(T))_{\rm red}$ and $Y+Z=\of^*(T)$, the sections do not
hit the divisor $Z$. We consider the blow up $a: \oX'\to \oX$ of
the sections $\overline{S_\alpha}$. Hence $\oX'$ is non-singular, the
divisor $Y'=a^{-1}Y$ is again a normal crossings divisor,
and $(\of\circ a)^*T$ is reduced in a neighborhood of
the exceptional divisors
$$
E_\alpha = a^{-1}\overline{S_\alpha}.
$$
$\iota$ acts on the relative Zariski tangent space of a section
$\overline{S_\alpha}$ by multiplication with $-1$, hence it
induces an action on $\oX'$, again denoted by $\iota$. The
restriction of $\iota$ to $E_\alpha$ is trivial, and $\iota$ acts
fixed point free on $\oX'\setminus \cup E_\alpha$. One has
$$R^i (\of\circ a)_*(\Omega^\bullet_{\oX'/\oS}(\log a^{-1}Y))
= \sH^i\oplus \sT_i, $$
where $\sT_i$ is an algebraically trivial bundle on which
$\iota$ acts trivially, and $\sT_i=0$ for $i$ odd.

The quotient, $\mathfrak{K}=\oX/(\iota)$ is non-singular and $h:
\mathfrak{K}\to \oS$ is a proper family, smooth over $S\subset \oS$.
Here $\mathfrak{K}$ stands for Kummer.

The ramification locus $\cup E_\alpha$ of $\oX'\to \mathfrak{K}$
is contained in the smooth locus of $\oX'\to \oS$, hence
$h^{-1}(T)$ is a normal crossings divisor, reduced in a
neighborhood of the image of $\cup E_\alpha$. So
for $\kappa=(h^{-1}(T))_{\rm red}$ the sheaf
$\Omega^1_{\mathfrak{K}/\oS}(\log \kappa)$ remains
locally free and
$$c_g(\Omega^1_{\mathfrak{K}/\oS}(\log \kappa))=
c_g(h^*\of_*(\Omega^1_{\oX/\oS}(\log Y)))
$$
in a neighborhood  of the multiple locus of $h^{-1}(T)$. Altogether,
$h:\mathfrak{K}\to \oS$ satisfies the assumption 1. - 3. in
Theorem \ref{GRR}, except possibly that the components of
$\kappa$ might be singular. As at the very end of the proof of
Theorem \ref{thmFC}, we blow up non-singular strata of 
$\kappa_{{\rm red}}$ to have 1. as well.

The Gau{\ss}-Manin bundles $\sH^i_\mathfrak{K}:=
R^ih_*\Omega^\bullet_{\mathfrak{K}/\oS} (\log \kappa)$
vanish for $i=2p+1, \ p>0$, and fulfill
\begin{gather}
\sH^{2p}_\mathfrak{K}= (R^{2p}(\of \circ a)_*
\Omega^\bullet_{\oX'/\oS}(\log a^{-1}Y))^{\iota},\notag
\end{gather}
where $^\iota$ means the invariants under $\iota$. Since $\iota$
acts trivially on $\overline{\sH}^{2p}$,
one obtains
$$
\sH^{2p}_\mathfrak{K}= \overline{\sH}^{2p} \oplus \sT_{2p}.
$$
Corollary \ref{char0} implies
\begin{gather}
{\rm ch}(\sum_{p\ge 0} \sH^{2p}_\mathfrak{K})=
(={\rm ch}(\sum_{p\ge 0} \overline{\sH}^{2p}))
\in CH^0(\oS)\otimes \Q.  \notag
\end{gather}
Lemma \ref{calcul} implies then
\begin{gather}
{\rm ch} (\overline{\sH}^1) \in CH^0(\oS)\otimes \Q. \notag
\end{gather}
This concludes the proof.
\end{proof}

\begin{thm} \label{charp}
Let $k$ be a field of characteristic $p \neq 2$.
Let $\of: \oX\to
\oS$ be a compactified family of abelian varieties as in Theorem
\ref{thmFC}. Then one has
\begin{gather}
{\rm ch}(\E^\vee \oplus \E)=
{\rm ch} (R^1\of_*\sO_\oX \oplus \of_*\Omega^1_{\oX/\oS}(\log
Y)) \in CH^0(\oS)\otimes \Q, \notag \end{gather}
or equivalently
$$ {\rm ch}_{2\ell} (\E)= {\rm ch}_{2\ell} (\E^\vee) =0 
\in CH^{2\ell}(\oS) \otimes \Q, \ \text{for} \ \ell\ge 1.
$$
\end{thm}
\begin{proof}
We replace in the proof of Theorem \ref{mainthm} the
Gau{\ss}-Manin bundle $\sH^i$ by the sum
$\sum_{j}R^j\of_*\Omega^{i-j}_{\oX/\oS}( \log Y)$ of the Hodge
bundles (see Theorem \ref{thmFC}, 6.).
\end{proof}
\bibliographystyle{plain}

\begin{thebibliography}{99}
\bibitem{Ba} Bateman, H. (compiled by the Bateman Manuscript
Project): Higher transcendental functions, vol.
I, McGraw-Hill Book Company (1953).
\bibitem{BE} Bloch, S.,  Esnault, H.:
Algebraic Chern-Simons theory.  Amer. J. Math.  {\bf 119}  (1997)
903--952.
\bibitem{BS} Borel, A., Serre, J.-P.: Le th\'eor\`eme de Riemann-Roch
(d'apr\`es Grothendieck), Bull. Soc. Math. France {\bf 86} (1958), 97--136.
\bibitem{D} Deligne, P.: \'Equations diff\'erentielles \`a points singuliers
r\'eguliers. Lecture Notes in Mathematics, {\bf 163}
Springer-Verlag, Berlin-New York, 1970. iii+133 pp.
\bibitem{D2} Deligne, P.: Th\'eorie de Hodge II, Publ. Math.
Inst. Hautes \'Et. Sci. {\bf 40} (1972) 5--57.
\bibitem{EV} Esnault, H.  Viehweg, E.: Logarithmic de~Rham
complexes and vanishing theorems.  Invent. Math.
{\bf  86 } (1986)
161--194.
\bibitem{E} Esnault, H.: Recent developments on characteristic
classes of flat bundles on
complex algebraic manifolds. Jahresber. Deutsch.
Math.-Verein.
{\bf 98 } (1996)
182--191.
\bibitem{FC} Faltings, G., Chai, C.-L.: Degeneration of Abelian
Varieties, Ergebnisse der Mathematik und ihrer Grenzgebiete,
3-te Folge,
{\bf 22} (1990), Springer Verlag.
\bibitem{Gr} Griffiths, P.: Periods of integrals on algebraic manifolds.
III. Some global differential-geometric properties of the period mapping.
 Inst. Hautes \'Etudes Sci. Publ. Math. {\bf  38} ( 1970) 125--180.
\bibitem{Gri} Griffiths, P. (Editor): Topics in transcendental algebraic
geometry. Ann of Math. Stud. {\bf 106}, Princeton Univ. Press,
Princeton, NJ. (1984)
\bibitem{H} Hirzebruch, F.: Topological Methods in Algebraic
Geometry, Grundlehren {\bf 131}, Springer Verlag, 3-rd edition (1966).
\bibitem{I} Iyer, J.: preprint Essen (2001).
\bibitem{GIT} Mumford, D.: Geometric Invariant Theory,
Ergebnisse der Mathematik und ihrer Grenzgebiete, 2-te Folge,
{\bf 34} (1982), Springer Verlag.
\bibitem{Mu} Mumford, D.: Towards an enumerative geometry of
the moduli space of curves.
Arithmetic and geometry, Vol. II,  271--328, Progr. Math.,{\bf  36},
Birkh\"auser Boston, Boston, MA, 1983.
\bibitem{vdG} van der Geer, G.: Cycles on the moduli space of
abelian varieties.  Moduli of curves and abelian varieties,  65--89,
Aspects Math., {\bf E33} , Vieweg, Braunschweig, 1999.
\bibitem{vdG2} van der Geer, G.: The Chow ring of the moduli
 space of abelian threefolds.  J. Algebraic Geom.
{\bf 7}  (1998) 753--770.
\end{thebibliography}
\renewcommand\refname{References}

\end{document}